\input amstex
\documentstyle{amsppt}
\magnification=1200 \hsize=13.8cm \catcode`\@=11
\def\NoLogo{\let\logo@\empty}
\catcode`\@=\active \NoLogo

\def\heat{\lf(\frac{\p}{\p t}-\Delta\ri)}

\def\lf{\left}
\def\ri{\right}

\def\e{\epsilon}
\def\p{\partial}

\def\R{\Bbb R}

\def\vp{\varphi}

\def \D {\Delta}

\documentstyle{amsppt}
\magnification=1200 \hsize=13.8cm \vsize=19 cm

\leftheadtext{Lei Ni} \rightheadtext{Entropy formula for linear
heat equation} \topmatter
\title{The entropy formula for linear heat equation}\endtitle

\author{Lei Ni\footnotemark }\endauthor
\footnotetext"$^{1}$"{{\it Math Subject Classifications.} Primary
58G11.} \footnotetext"$^{2}$"{{\it Key Words and Phrases.} Entropy
formula, heat equation, logarithmic Sobolev.}
\footnotetext"$^{3}$"{Research partially supported by NSF grant
DMS-0328624, USA.}
\address
Department of Mathematics, University of California, San Diego, La
Jolla, CA 92093
\endaddress
\email{ lni\@math.ucsd.edu}
\endemail

\affil { University of California, San Diego}
\endaffil

\date  June 2003 (Revised August 2003)\enddate

\abstract We derive the entropy formula for the linear heat
equaiton on general Riemannian manifolds and prove that it is
monotone non-increasing on manifolds with nonnegative Ricci
curvature. As applications, we study the relation between the
value of entropy and the volume of balls of various scales. The
results are simpler version, without Ricci flow, of Perelman's
recent results on volume non-collapsing for Ricci flow on compact
manifolds. We also prove that if the entropy achieves the maximum
value zero on any complete Riemannian manifold with nonnegative
Ricci curvature, then the manifold must be isometric to the
Euclidean space.

\endabstract

\endtopmatter

\document

\subheading{\S0 Introduction}\vskip .2cm

In a recent paper of Perelman\cite{P1}, an entropy formula for
Ricci flow was derived. The formula turns out being of fundamental
importance in the study of Ricci flow (cf. \cite{P1, Sections 3,
4, 10}) as well as the K\"ahler-Ricci flow \cite{P2}. The
derivation of the entropy formula in \cite{P1, Section 9}
resembles the gradient estimate for the linear heat equation
proved by Li-Yau in another fundamental paper \cite{L-Y} on the
linear parabolic equation. This suggests that there may exist a
similar entropy formula for the linear heat equation. The purpose
of this short note is to show  such entropy formula  and derive
some applications of the new entropy formula.

Let $M$ be a complete Riemannian manifold. We study the heat
equation
$$
\heat u =0. \tag 0.1
$$
Following \cite{P1}, we define
$$
{\Cal W}(f, \tau)=\int_M \lf(\tau |\nabla
f|^2+f-n\ri)\frac{e^{-f}}{(4\pi \tau)^{\frac{n}{2}}}\, dv \tag 0.2
$$
restricted to $(f, \tau)$ satisfying
$$
\int_M\frac{e^{-f}}{(4\pi \tau)^{\frac{n}{2}}}\, dv=1 \tag 0.3
$$
with $\tau >0$.

\proclaim{Theorem 0.1} Let $M$ be a closed Riemannian manifold.
Assume that $u$ is a positive solution to the heat equation (0.1)
with $\int_M u\, dv =1$. Let $f$ be defined as
$u=\frac{e^{-f}}{(4\pi \tau)^{\frac{n}{2}}}$ and $\tau=\tau(t)$
with $\frac{d \tau}{d t} =1$. Then
$$
\frac{d {\Cal W}}{d t}= -\int_M 2\tau \lf(|\nabla_i\nabla_j
 f-\frac{1}{2\tau}g_{ij}|^2+R_{ij}f_{i}f_{j}\ri)u\, dv. \tag 0.4
$$
In particular, if $M$ has nonnegative Ricci curvature, ${\Cal
W}(f,\tau)$ is monotone decreasing along the heat equation.
\endproclaim

Notice that in the case that $M$ is Ricci flat, the result above
is indeed a special case of Perelman's result. We just show that
it in fact also holds for all metrics with nonnegative Ricci
curvature.

The result can be derived out of a point-wise differential
inequality. The proof of Theorem 0.1 and the argument of \cite{P1}
gives the following differential inequality for the fundamental
solution to the heat equation.

\proclaim{Theorem 0.2} Let $M$ be a closed Riemannian manifolds
with nonnegative Ricci curvature. Let $H$ be the positive heat
kernel. Then
$$
t(2\D f-|\nabla f|^2)+f-n \le 0, \tag 0.5
$$
for $t>0$ with $H=\frac{e^{-f}}{(4\pi t)^{\frac{n}{2}}}$.
\endproclaim

Notice that Li-Yau's gradient estimate
$\frac{u_t}{u}-\frac{|\nabla u|^2}{u^2}+\frac{n}{2t}\ge 0$ is
equivalent to
$$
t(2\D f)-n\le 0. \tag 0.6
$$
The inequality (0.6) can be viewed as a generalized Laplacian
comparison theorem. In deed,  the Laplacian comparison theorem on
$M$  is a consequence of (0.6) by applying the inequality to the
heat kernel and letting $t\to 0$. This suggestions that one can
view $\bar{L}(x,t)=4tf$ as the square of a time-dependent
`distance function'. Then (0.6), which says that $\Delta
\bar{L}\le 2n $,  simply generalizes the standard Laplacian
comparison $\D r^2\le 2n$ on any Ricci non-negative manifold to
such generalized `distance function'. From this point of view, one
can think (0.5) as a Laplacian comparison theorem in the
space-time since it says $\D \bar{L}+\bar{L}_t\le 2n$. The similar
inequality \cite{P1, (7.15)} was one of the important new
discoveries of Perelman. Applying similar consideration, one can
also view the entropy estimate in \cite{P1} as a generalization of
the space-time Laplacian comparison theorem. This is also related
to the reduced volume monotonicity of \cite{P1}. It was pointed to
us later by Professor S. T. Yau that Hamilton and him also noticed
(0.5) for the Ricci flow a few years ago. It seems that  (0.5) and
(0.6) do not imply each other. It is interesting to find out if
there is any deeper connection between them.

For closed manifolds, following  \cite{P1},  one can define
$$
\mu(\tau)=\inf_{\int_M u\, dv=1}{\Cal W}(f, \tau). \tag 0.7
$$
A direct consequence of Theorem 0.1 is the following

\proclaim{Corollary 0.1} On manifolds with nonnegative Ricci
curvature, $\mu(\tau)$ is a monotone decreasing function of
$\tau$. Moreover $\mu(\tau)\le 0$ with $\lim_{\tau\to
0}\mu(\tau)=0$.
\endproclaim
Whenever it makes sense (for example, when $M$ is
simply-connected,  negatively curved with lower bound on its
curvature), as in \cite{P 1} one can also define
$\nu=\inf_{\tau}\mu(\tau)$. It can be thought as some sort of
isoperimetric constant. When $M=\R^n$, $\nu=0$. Thanks to the
gradient estimates of Li-Yau \cite{L-Y}, the above results still
hold on complete noncompact manifolds with nonnegative Ricci
curvature. As an application of the entropy formula obtained in
Theorem 0.1  we prove the following result.

\proclaim{Theorem 0.3} Let $M$ be a complete Riemannian manifold
with nonnegative Ricci curvature. Then $\mu(\tau)\ge 0$ for some
$\tau>0$ if and only if $M$ is isometric to $\R^n$.
\endproclaim

In \cite{G} (see also \cite{S, W}), a sharp logarithmic Sobolev
inequality (in different disguises in \cite{S, W}) was proved on
$\R^n$. When $M=\R^n$, the inequality is equivalent to
$$
\int_M \lf(\frac{1}{2}|\nabla
f|^2+f-n\ri)\frac{e^{-f}}{(2\pi)^{\frac{n}{2}}}\, dv \ge 0 \tag
0.8
$$
for all $f$ with $\int_M\frac{e^{-f}}{(2\pi)^{\frac{n}{2}}}\,
dv=1$.

Since (0.8) is equivalent to $\mu(\frac{1}{2})\ge 0$, a simple
corollary of Theorem 0.3 is the following result on the relation
between the logarithmic Sobolev inequality and the geometry of the
manifolds, which is originally due to Bakry, Concordet and Ledoux
\cite{B-C-L}.

\proclaim{Corollary 0.2} Let $M$ be a complete Riemannian manifold
with nonnegative Ricci curvature. Then (0.8) holds on $M$ if and
only if $M$ is isometric to $\R^n$.
\endproclaim

It can be shown  that the (0.8) holds on any manifold with sharp
isoperimetric inequality, or equivalently the sharp $L^1$-Sobolev
inequality. Under the request of some readers of the preliminary
version, we include a proof of this fact (cf. Proposition 3.1)
using the spherical symmetrization. The proof was communicated to
the author by Perelman last year in November. The proof only uses
the
 spherical symmetrization to compare with the Euclidean case. One can find the
  simple elegant proof of (0.8) by
 Beckner and Pearson in
 \cite{B-P}, which makes use of the fact that the product of
 Euclidean spaces is still Euclidean together with the sharp $L^2$-Sobolev inequality.

It turns out that ${\Cal W}(f, \tau)$ being finite, where
$u=\frac{e^{-f}}{(4\pi t)^{\frac{n}{2}}}$ is the heat kernel, also
has strong geometric and topological consequences. For example, in
the case $M$ has nonnegative Ricci curvature, it implies that $M$
has finite fundamental group. In fact we can show that

\medskip

{\it $M$ is of maximum volume growth if and only if the entropy
${\Cal W}(f, t)$ is uniformly bounded for all $t>0$, where
$u=\frac{e^{-f}}{(4\pi t)^{\frac{n}{2}}}$ is the heat kernel.}

\medskip

The  analogy of above  was discovered originally in \cite{P1} for
the ancient solution to Ricci flow with bounded nonnegative
curvature operator, where he claimed that an ancient solution to
the Ricci flow with nonnegative curvature operator is
$\kappa$-non-collapsed if and only if the entropy is uniformly
bounded for any fundamental solution to the conjugate heat
equation.

Without assuming the nonnegativity of the Ricci curvature, the
bound on $\mu(\tau)$ also implies the uniform lower bound on the
volume of balls of certain scales. Namely, it implies the volume
noncollapsing, as in the $\kappa$-noncollapsing theorem of
Perelman \cite{Theorem 4.1, P1}, therefore an uniform upper bound
of the diameter, if the manifold has finite volume. In some sense,
$\mu(\tau)$ reflects the isoperimetric property of $M$ for the
scale parametrized by $\tau$.

\medskip

{\it Acknowledgement.} The author would like to thank Professor
Ben Chow for encouragement to study Perelman's recent papers and
discussions,  Professor Peter Li for helpful suggestions. Special
thanks goes to Professor Perelman since most results in this paper
are the simpler versions  of the corresponding results in
 \cite{P1} without Ricci flow.

 \subheading{\S 1 Proof of Theorem 0.1 and 0.2}

\medskip

We start with the following two lemmas.

\proclaim{Lemma 1.1} Let $M$ be a complete Riemannian manifold.
Let $u$ be a positive solution to (0.1). Then
$$
\heat w = -2 \lf(|\nabla_i \nabla_j f|^2+R_{ij}f_if_j\ri)
-2<\nabla w, \nabla f> \tag 1.1
$$
where $w=2\D f -|\nabla f|^2$ and $f=-\log u$.
\endproclaim
\demo{Proof} Direct calculation shows that
$$
\split \heat w&= -2|\nabla_i \nabla_j f|^2-2R_{ij}f_i f_j-2<\nabla
(\D f), \nabla f>+2<\nabla (f_t),\nabla f>\\
&\ \ -2\D(f_{t})+2f_{tt} \\
&= -2|\nabla_i \nabla_j f|^2-2R_{ij}f_i f_j-2<\nabla(|\nabla
f|^2), \nabla f>-2(|\nabla f|^2)_t\\
& = -2|\nabla_i \nabla_j f|^2-2R_{ij}f_i f_j-2<\nabla(|\nabla
f|^2+2f_t), \nabla f>\\
&=-2|\nabla_i \nabla_j f|^2-2R_{ij}f_i f_j-2<\nabla w, \nabla f>.
\endsplit
$$
Here we have used $w=2f_t+|\nabla f|^2$ and $\heat f=-|\nabla
f|^2$.
\enddemo
\proclaim{Lemma 1.2} Let $M$ and $u$ be as in Lemma 1.1. Let $ W =
\tau(2\D f-|\nabla f|^2)+f-n$, where we write
$u=\frac{e^{-f}}{(4\pi \tau)^{\frac{n}{2}}}$. Here $\tau=\tau(t)$
with $\frac{d \tau}{dt}=1$. Then
$$
\heat W = -2\tau\lf(|\nabla _i\nabla_j
f-\frac{1}{2t}g_{ij}|^2+R_{ij}f_if_j\ri)-2<\nabla W, \nabla f>.
\tag 1.2
$$
\endproclaim
\demo{Proof} One can proceed directly. Here we use Lemma 1.1 to
simplify the calculation a little. Let $\bar{f}=-\log u$. Then we
have that
$$
W=\tau w+\bar{f}-\frac{n}{2}\log (4\pi \tau)-n.
$$
Keep in mind  $\heat \bar f =-|\nabla \bar f|^2$. The direct
calculation shows that
$$
\split \heat W &= \tau \heat w+w-|\nabla \bar f|^2-\frac{n}{2 \tau}\\
&= -2\tau|\nabla_i\nabla_j f|^2-2\tau R_{ij}f_if_j-2\tau<\nabla w,
\nabla f>+|\nabla \bar
f|^2+2{\bar f}_t-|\nabla \bar f|^2-\frac{n}{2\tau}\\
& = -2\tau|\nabla_i\nabla_j f|^2-2<\nabla W, \nabla f>+2\D \bar f
-\frac{n}{2\tau}-2\tau R_{ij}f_if_j\\
&= -2\tau|\nabla_i\nabla_j f-\frac{1}{2\tau}g_{ij}|^2-2<\nabla W,
\nabla f>-2\tau R_{ij}f_if_j.
\endsplit
$$
Here we have used $\nabla f =\nabla \bar f$.

\proclaim{Remarks} 1. Lemma 1.1 has its corresponding version for
the Ricci flow. Namely, if $g_{ij}$ satisfies the back-ward Ricci
flow equation $\frac{\p}{\p t}g_{ij}=2R_{ij}$ and $u$ is a
solution to $(\frac{\p}{\p t}-\D +R)u=0$. Define $w=(2\D f-|\nabla
f|^2+R)$. Then
$$
\heat w=-2|R_{ij}+f_{ij}|^2-2<\nabla w, \nabla f>. \tag 1.3
$$
Here $u=e^{-f}$. One can easily see that (1.3) implies the formula
(1.2) of \cite{P1}. This also gives another derivation of the
first monotonicity formula in \cite{P1}.

2. The above approach of the proof to Lemma 1.2 was motivated by
the statistical analogy in \cite{P1, Section 6}. One can also use
the similar approach to derive Proposition 9.1 of \cite{P1} from
(1.3) above. This would simplify the calculation a little and
reflect the relation between the energy and the entropy quantity.
\endproclaim

\enddemo

\demo{Proof of Theorem 0.1} The proof of Theorem 0.1 follows from
the simple observation $u \nabla f=-\nabla u$, therefore
$$
\split
 \heat \lf(Wu\ri)& =-2\tau\lf(|\nabla _i\nabla_j
f-\frac{1}{2\tau}g_{ij}|^2+R_{ij}f_if_j\ri)u\\
&\ \ -2<\nabla W, \nabla f>u-2<\nabla W, \nabla u>\\
&=-2\tau\lf(|\nabla _i\nabla_j
f-\frac{1}{2\tau}g_{ij}|^2+R_{ij}f_if_j\ri)u,
\endsplit
$$
and integration by parts.
\enddemo

\demo{Proof of Theorem 0.2} We can apply Perelman's argument in
the proof of Corollary 9.3 of \cite{P1}. For any $t_0>0$, let $h$
be any positive function. We solve the backward heat equation
starting from $t_0$ with initial data $h$. Then we have that
$$
\split \frac{d}{d t}\int_M hWu\, dv &= \int_M (h_t)(Wu)+h(Wu)_t\,
dv \\
&= \int_M \lf(h_t+\Delta h\ri)Wu +h((Wu)_t-\D (Wu))\, dv\\
&\le 0.
\endsplit
$$
Using the fact that $(\int_M hWu\, dv)|_{t=0}=0$, when $u$ is the
fundamental solution,  we have that
$$
\int_M h(Wu)\, dv \le 0
$$
for any $t_0>0$ and any positive function $h$. This implies that
$Wu\le 0$. Therefore $W\le 0$.
\enddemo

\subheading{\S 2 Extensions and the value of $\mu(0)$}

\medskip

The first extension is to complete noncompact manifolds. From the
proof, it is easy to see that Theorem 0.1 and 0.2 hold as long as
the integration by parts can be justified. We focus on the case
$M$ has nonnegative Ricci curvature. Since we have the gradient
estimate of Li-Yau for the positive solutions one indeed can make
the integration by part rigorous, keeping in mind that $u$ is
assumed integrable in our consideration of the entropy. One of the
reference where one can find the estimates on derivatives of $u$
is \cite{C-N, Section 3}.

Another extension of Theorem 0.1 and Theorem 0.2 is for manifolds
with boundary. In this case, it is not hard to show that the
theorem holds on manifolds with convex boundary. In fact, in this
case
$$
\split
 \frac{\p W}{\p \nu}& =(2\tau f_{\tau} +\tau |\nabla
f|^2+f-n)_{\nu}\\
&= 2\tau \sum_i^{n-1}f_{i}f_{i\nu}\\
&=-2\tau h_{ij}f_if_j\\
&\le 0.
\endsplit
$$
Here $h_{ij}$ denotes the second fundamental form of $\p M$.
Therefore, Theorem 0.1 and Theorem 0.2 hold for positive solution
 $u$ with the Neumann boundary condition $\frac{\p u}{\p \nu}=0$.

 \proclaim{Corollary 2.1} Let $M$ be a compact manifold with
 boundary. Let $u$ be a positive solution to (0.1) with the Neumann boundary condition.
 Let $f$ and $\tau$ be as in Theorem 0.1. Then
 $$
\frac{d}{d t}{\Cal W}=-\int_M 2\tau \lf(|\nabla_i\nabla_j
 f-\frac{1}{2\tau}g_{ij}|^2+R_{ij}f_{i}f_{j}\ri)u\, dv-2\int_{\p M}\tau II((\nabla f)^{T},
 (\nabla f)^{T})\,
 dA.
 \tag 2.1
 $$
 Here $II(\cdot,\cdot )$ is the second fundamental form of $\p M$
 and $(\nabla f)^{T}$ is the tangential projection of $\nabla f$
 on $\p M$.
In particular, in the case $M$ has nonnegative Ricci and $\p M$ is
convex, ${\Cal W}$ is monotone decreasing. Moreover, if $u$ is the
fundamental solution,
$$
t(2\D f-|\nabla f|^2)+f-n \le 0, \tag 2.2
$$
for $t>0$.
 \endproclaim

Since we know that $\mu(\tau)$ is monotone, it is nice to know the
value of $\mu(\tau)$ as $\tau \to 0$. We can adapt the argument of
\cite{P1} to prove that $\lim_{\tau \to 0}\mu(\tau)=0$. Since the
argument in \cite{P1} of this part is  sketchy we include a
detailed proof for the sake of the reader. But the original idea
is certainly from \cite{P1}.

\proclaim{Proposition 2.1} Let $M$ be a closed manifold.
$$
\lim_{\tau \to 0}\mu(\tau)=0. \tag 2.3
$$
\endproclaim
\demo{Proof} It is easy to see that $\mu(\tau)\le 0$ by Theorem
0.2. Assume that there exists $\tau_k\to 0$ such that
$\mu(\tau_k)\le c<0$ for all $k$. We show that this will
contradict the logarithmic Sobolev inequality. We are going to
blow up the metric by $\frac{1}{2}\tau^{-1}$. First we can
decompose $M$ into open subsets $U_1,U_2, \cdots, U_N$ such that
each $U_j$ is contained inside some normal coordinates and each
$U_j$ also contains $B(o_j, \delta)$, a ball of radius $\delta$,
for some small $\delta>0$. Now let
$g^\tau=\frac{1}{2}\tau^{-1}g_{ij}$ and $g_k=g^{\tau_k}$. It is
clear that $(U_j, g_k, o_j)$ converges to $(\R^n, g_0, 0)$ in
$C^{\infty}$ norm. We will also identify the compact subset of
$\R^n$ with the compact subset of $U_j$.

It is easy to see that
$$
{\Cal W}_{g}(f, \tau)=\int_M (\frac{1}{2}|\nabla
f|^2_{\tau}+f-n)\frac{e^{-f}}{(2\pi)^{\frac{n}{2}}}\, dv_\tau
$$
with restriction $\int_M \frac{e^{-f}}{(2\pi)^{\frac{n}{2}}}\,
dv_\tau=1$, where $|\cdot|_{\tau}$ is the norm with respect to
$g^{\tau}=\frac{1}{2}\tau^{-1}g$ and $dv_{\tau}$ is the
corresponding volume form. It is also convenient to write in more
standard form:
$$
{\Cal W}(\psi, \tau)=\int_M \lf(2|\nabla \psi|^2_{\tau}-(\log
\psi^2)\psi^2-(\frac{n}{2}\log(2\pi)+n)\psi^2\ri)\, dv_{\tau} \tag
2.4
$$
restricted to $\int_M \psi^2\, dv_{\tau}=1$. Let $\vp_k$ be the
minimizer realizing $\mu(\tau_k)$. Then we have that
$$
-2\D_k \vp_k -2\vp_k \log \vp_k
=\lf(\mu(\tau_k)+n+\frac{n}{2}\log(2\pi)\ri)\vp_k \tag 2.5
$$
and
$$
\int_M \vp_k^2\, dv_k =1. \tag 2.6
$$
Here $\D_k$ denote the Laplacian of $g_{\tau_k}$ and $dv_k
=dv_{\tau_k}$. Due to the monotonicity, we can also assume that
$\mu(\tau_k)\ge -A$ for some $A>0$ independent of $k$.

Now we write $F_k(\psi)=2|\nabla \psi|^2_{\tau_k}-(\log
\psi^2)\psi^2-(\frac{n}{2}\log(2\pi)+n)\psi^2$. It is a easy
matter to check that
$$
\frac{\int F(\lambda \phi)\, dv_{\tau}}{\int (\lambda \psi)^2\,
dv_{\tau}}=\frac{\int F(\psi)\, dv_{\tau}}{\int \psi^2\,
dv_\tau}-\log \lambda^2. \tag 2.7
$$
By the assumption that $\mu(\tau_k)\le c<0$ we know that
$$
\int_M F(\vp_k)\, dv_{k} \le c<0.
$$
By passing to subsequence we can assume that
$$
\int_{U_1}F(\vp_k)\, dv_{k} \le \frac{c}{N}<0.
$$
It is easy to see that $\int_{U_1}\vp_k^2\, dv_k \le 1$. Combining
the above  with (2.5) and the fact that $g_{k}$ converges to $g_0$
on every fixed compact subset of $\R^n$, the elliptic PDE theory
implies  that there exists a subsequence of $\vp_k$ (still denote
by $\vp_k$) such that it converges uniformly on every compact
subset of $\R^n$. If the limit function $\vp_{\infty}$ exists and
$\int_{\R^n} \vp_{\infty}^2\, dv_0 > 0$ we claim that we will get
contradiction to the logarithmic Sobolev inequality (0.8). In fact
in this case we just denote $\e^2=\int_{\R^n}\vp^2_{\infty}$.
Clearly $0<\e\le 1$ by the assumption. Since
$$
\int_{\R^n} F(\vp_{\infty})\, dv_0\le \frac{c}{N}<0 \tag 2.8
$$
then by (2.7) we have that
$$
\int_{\R^n}F(\frac{1}{\e}\vp_{\infty})\, dv_0 \le
\frac{c}{N}+2\log \e <\frac{c}{N}.
$$
Let
$\frac{e^{-f_{\infty}}}{(2\pi)^{\frac{n}{2}}}=\lf(\frac{1}{\e}\vp_{\infty}\ri)^2$
we have that
$\int_{\R^n}\frac{e^{-f_{\infty}}}{(2\pi)^{\frac{n}{2}}}\, dv_0=1$
and
$$
\int_{\R^n}(\frac{1}{2}|\nabla
f_{\infty}|^2_{0}+f_{\infty}-n)\frac{e^{-f_{\infty}}}{(2\pi)^{\frac{n}{2}}}\,
dv_0 \le \frac{c}{N}<0.
$$
This is a contradiction to the sharp  logarithmic Sobolev
inequality (0.8). On the other hand, if $\e=0$.  which would imply
$\vp_{\infty}=0$. This contradicts (2.8). We therefore  complete
the proof of the proposition.

\enddemo

\subheading{\S 3 Bounded entropy and volume growth}

\medskip

The main purpose of the section is to prove Theorem 0.3 and show
that for the manifold with nonnegative Ricci curvature the
finiteness of the entropy for the heat kernel is equivalent to the
manifold has maximum volume growth. We first include a short
discussion on the logarithmic Sobolev inequality. We say $M$ has
logarithmic Sobolev inequality if
$$
\int_M \lf(\frac{1}{2}|\nabla
f|^2+f\ri)\frac{e^{-f}}{(2\pi)^{\frac{n}{2}}}\, dv\ge -C_1 \tag
3.1
$$
for all $f$ with restriction $\int_M
\frac{e^{-f}}{(2\pi)^{\frac{n}{2}}}\, dv=1$. This is equivalent to
the finiteness of $\mu(\frac{1}{2})$. It is an easy matter to see
that the regular Sobolev inequality implies (3.1).
 In particular, it holds on minimal submanifolds in
$\R^n$, which is a special case of the general result in \cite{E}.
Since one has the $L^2$-Sobolev inequality on a closed manifold,
 any closed manifold satisfies (3.1).
In this case the dependence of the constant $C_1$ can be
explicitly traced, applying Lemma 2 of \cite{L1}. This would  in
turn gives the explicit dependence  of the $\kappa$ constant on
the geometry of the initial metric in the $\kappa$ noncollapsing
theorem of \cite{P1, Section 4}.

The inequality (3.1) is also equivalent to the ultracontractivity
as pointed out in \cite{D}, which then follows from conditions
such as the lower bound on the Ricci curvature (i.e. $Ric\ge -K$
for some $K\ge 0$) and $\inf V_x(1)\ge \delta>0$. Therefore,
$\mu(\frac{1}{2})$ is finite for a large class of manifolds.

It is interesting to find out on which manifolds the logarithmic
Sobolev inequality holds with $C_1=n$. Namely (0.8) holds. It was
pointed out in \cite{P1} that the sharp isoperemetric inequality
also implies the sharp logarithmic Sobolev inequality. The
following was the proof suggested by the  communication with
Perelman. The result was somewhat conjectured in \cite{R}.

\proclaim{Proposition 3.1 (Perelman)} Let $M$ be a complete
manifold such that
$$
A(\p \Omega)\ge c_n V(\Omega)^{\frac{n-1}{n}},
$$
for any compact domain $\Omega$ with the Euclidean constant $c_n$.
Here $A(\p \Omega)$ is the area of the boundary $\p \Omega$. Then
(0.8) holds on $M$
\endproclaim
\demo{Proof} As we know in the proof of Proposition 2.1, (0.8) is
equivalent to
$$
\int_M 2|\nabla \phi|^2-(\log\phi^2)\phi^2-(\frac{n}{2}\log
(2\pi)+n)\phi^2\, dv \ge 0. \tag 3.2
$$
It suffices to prove the result for compact supported nonnegative
function $\phi$. Let $M'=\{x|\phi(x)>0\}$. Let $\bar{B}(R)$ be a
ball of radius $R$ in $\R^n$ such that $Vol(\bar{B}(R))=Vol(M')$.
We define that $F(t)=Vol(\{x\in M'| \phi(x)\ge t\})$. We also
denote $M_t=\{x\in M'|\phi(x)\ge t\}$. $\Gamma_t =\p M_t$. Let
$g(|y|)$ be a function on $\bar{B}(R)$ such that $Vol(\{y|g(y)\ge
t\})=F(t)$ and $g(R)=0$. We can define $\bar{M}_t$ and
$\bar{\Gamma}_t$ similarly. Clearly $Vol(M_t)=Vol(\bar{M}_t)$ and
$A(\Gamma_t)\ge A(\bar{\Gamma}_t)$ by the isoperimetric
inequality. The simply integration by parts shows that
$$
\int_0^\infty \lambda'(s) F(s)\, ds =\int_{M'}\lambda(f)\, dv \tag
3.3
$$
for any Lipschitz function $\lambda (t)$ with  $\lambda (0)=0$.
This implies that
$$
\int_{M'} (\log \phi^2)\phi^2+(\frac{n}{2}\log(2\pi)+n)\phi^2\, dv
=\int_{\bar{B}(R)}(\log g^2)g^2+(\frac{n}{2}\log (2\pi)+n)g^2 \,
d\bar{v}.
$$
On the other hand the isoperimetric inequality implies
$$
\int_{M'}|\nabla \phi|^2\, dv \ge \int_{\bar{B}(R)}|\bar{\nabla}
g|^2\, d\bar{v}. \tag 3.4
$$
In fact, the co-area formula shows that
$$
\int_{M'}|\nabla \phi|^2\, dv=\int_0^\infty \int_{\Gamma_t}|\nabla
\phi|\, dA\, dt.
$$
and
$$
F(t)=\int_t^{\infty} \int_{\phi =s}\frac{1}{|\nabla \phi|}\, dA\,
ds.
$$
Combining with the fact $F(t)=Vol(\bar{M}_t)$ we have that
$$
\int_{\Gamma_t} \frac{1}{|\nabla \phi|}\,
dA=\int_{\bar{\Gamma}_t}\frac{1}{|\bar{\nabla} g|}\, d\bar{A}.
$$
Using the fact that $|\bar{\nabla} g|$ is constant on
$\bar{\Gamma}_t$, by H\"older inequality,  we have that
$$
\split
 \lf(\int_{\bar{\Gamma}_t}|\bar{\nabla} g|\, d\bar{A} \ri) \lf(\int_{\bar{\Gamma}_t}\frac{1}{|\bar{\nabla} g|}\, d\bar{A} \ri)
& =A^2(\bar{\Gamma}_t) \\
&\le A^2(\Gamma_t)\\
&\le \lf(\int_{\Gamma_t}|\nabla \phi|\, dA \ri)
\lf(\int_{\Gamma_t}\frac{1}{|\nabla \phi|}\, dA \ri)
\endsplit
$$
which then implies that
$$
\int_{\Gamma_t}|\nabla \phi|\, dA
\ge\int_{\bar{\Gamma}_t}|\bar{\nabla} g|\, d\bar{A}. \tag 3.5
$$
Since (3.4) implies (3.5) we completes the proof.
\enddemo

We should point out that in \cite{B}, Beckner provides a proof of
(0.8) from the isoperimetric inequality using the product
structure of the Euclidean spaces. The above argument  just
reduces the (0.8) for any manifolds, with the sharp isoperimetric
inequality, to the Euclidean space (with the same dimension) case.
It does not prove the Euclidean case itself. Now we prove Theorem
0.3.

 \demo{Proof of Theorem 0.3} By the assumption, one can find
 $\tau_0$ such that
$\mu(\tau_0)\ge 0$. But on the other hand, Theorem 0.2 implies
that ${\Cal W}(f,t)\le 0$ for $H=\frac{e^{-f}}{(4\pi
t)^{\frac{n}{2}}}$ being the heat kernel. This implies $\mu(t)=0$
for $0<t<\tau_0$.  Applying the equality case in Theorem 0.1 we
have that $f_{ij}-\frac{1}{2t}g_{ij}=0$, which implies that
$$
2t\D f = n. \tag 3.6
$$
On the other hand, by \cite{C-L-Y, V} we know that $\lim_{t\to
0}-4t\log H=r^2(x,y)$. In particular,
$$
\lim_{t\to 0} 4tf=r^2(x,y).
$$
Then (3.6) implies that
$$
\D r^2(x,y) =2n. \tag 3.7
$$
Combining with the assumption that Ricci is nonnegative this
implies that $M$ is isometric to $\R^n$. In fact from (3.7) one
can easily obtain that
$$
\frac{A_x(r)}{V_x(r)}=n
$$
where $A(r)$ and $V(r)$ denotes the area of $\p B_x(r)$ and the
volume of $B_x(r)$, which then implies that $V_x(r)$ is same as
the volume function of Euclidean balls. The equality case of the
volume comparison theorem implies $M=\R^n$

\enddemo

It is clear that Theorem 0.3 implies Corollary 0.2. The proof of
\cite{B-C-L} to Corollary 0.2 relies on a deep result of Peter Li
\cite{L2} on the large time behavior of the heat kernel. Since we
only uses the behavior of the heat kernel near $t=0$, our proof is
dual to theirs in some sense. In \cite{Le}, the author proved that
the sharp Sobolev inequalities on a Ricci nonnegative manifold
also implies the manifold is isometric to $\R^n$. The case of
$L^1$-Sobolev, which is equivalent to the isoperimetric inequality
in Proposition 3.1, is relatively simple. The other cases are more
involved. Please see \cite{Le} for details. It was also asked in
\cite{Le}  if the sharp Nash inequality implies the same
conclusion or not. That still remains open.

\proclaim{Proposition 3.2} Let $M$ be a complete Riemannian
manifold with nonnegative Ricci curvature. Assume that $M$ has
maximum volume growth, namely $\frac{V_o(r)}{r^n}\ge \theta$ for
some $\theta>0$. Then there exists $A=A(\theta, n)>0$ such that
$$
{\Cal W}(f, t)\ge -A \tag 3.8
$$
 for $u=\frac{e^{-f}}{(4\pi
t)^{\frac{n}{2}}}$ being the heat kernel. On the other hand, (3.8)
implies that $M$ has maximum volume growth. Namely
$\frac{V_o(r)}{r^n}\ge \theta$ holds for some $\theta=\theta(n,
A)$.
\endproclaim
\demo{Proof} Let $v=\sqrt{u}$. One can rewrite ${\Cal W}(f, t)$ as
$$
{\Cal W}=4t\int_M |\nabla v|^2\, dv -\int_M \log(v^2)v^2\, dv
-\lf(n+\frac{n}{2}\log (4\pi t)\ri). \tag 3.9
$$
On the other hand, by Li-Yau's heat kernel estimate
$$
\split
 v^2&=H(x,y, t)\le \frac{C(n)}{V_x(\sqrt{t})}\\
 &\le
\frac{C(n)\theta}{t^{\frac{n}{2}}}.
\endsplit
$$
Hence
$$
{\Cal W}\ge -\log(C(n)\theta)-n-\frac{n}{2}\log (4\pi).
$$
Here we have used the fact $\int_M v^2\, dv =1$.

To prove the second half of the claim we need to use the lower
bound estimate of Li-Yau as  well as the gradient estimate for the
heat kernel. We first estimate the first term in (3.9) using
inequality (0.6), the Li-Yau's gradient estimate.
$$
\split t\int_M |\nabla v|^2\, dv &= t\int_M \frac{|\nabla
H|^2}{H}\, dv\\
&\le t\int_M \lf(H_t+\frac{n}{2t}H\ri)\, dv\\
&=\frac{n}{2}.
\endsplit \tag 3.10
$$
The second term can be estimated as
$$
\split
 -\int_M \log (H)H \, dv&\le -\int_M
 \log\lf(\frac{C_5(n)}{V_x(\sqrt{t})}\exp(-\frac{r^2(x,y)}{3t})\ri)H\,
 dv_y\\
 &\le C_6(n)+\log(V_x(\sqrt{t})) +\frac{1}{3t}\int_Mr^2(x,y)H(x,y,
 t)\, dv_y\\
 &\le C_7(n) +\log(V_x(\sqrt{t})).
\endsplit
\tag 3.11
$$
Here $C_i$ are positive constants only depending on $n$. We also
have used Theorem 3.1 of \cite{N} to estimate the last term of the
second line above. Putting the assumption ${\Cal W}\ge -A$ and
(3.9)--(3.11) together we have the lower bound (3.1) for the
volume.
\enddemo

The similar result as above was claimed  in \cite{P1, Section 11}
for the Ricci flow ancient solutions.  The proof here is easier
than the nonlinear case considered in \cite{P1}. In fact,
Proposition 3.2 here  can be used in the proof of Theorem 10.1 of
\cite{P1}.

\subheading{\S 4 Manifolds with bounded $\mu(\tau)$}

\medskip

The following result gives the geometric implication of the
non-sharp logarithmic Sobolev inequality (3.3), or bounded
$\mu(\tau)$. The result can be thought as Riemannian version of
the $\kappa$ non-collapsing result of \cite{P1}. Notice that we do
not even require $M$ has nonnegative Ricci curvature. The results
in this section are in the line of  Perelman's work on
K\"ahler-Ricci flow \cite{P2}. However, the arguments in the
nonlinear case are technically more involved than the case treated
here, especially on the diameter bound.

\proclaim{Proposition 4.1} Let $M$ be a complete manifold. Assume
that $\mu(\tau)\ge -A$, for all $0\le \tau\le T$, for some
constant $A>0$. Then there exists a positive constant $\kappa(A,
n)>0$ such that
$$
V_x(R)\ge \kappa R^{n} \tag 4.1
$$
for all $R^2\le T$.  In particular,  if $\mu(\tau)\ge -A$ for all
$\tau\ge 0$, $M$ has at least the Euclidean volume growth.
\endproclaim
\demo{Proof} The first observation is that
$$
\mu(\tau)\le \int_M \tau 4|\nabla h|^2-\lf(\log
h^2+\frac{n}{2}\log (4\pi \tau)\ri)h^2\, dv \tag 4.2
$$
for compact supported nonnegative function $h$. If we have that
$$
V_x(\frac{R}{2})\ge \eta V_x(R) \tag 4.3
$$
for $\eta =\lf(\frac{1}{3}\ri)^{n}$ we will have the estimate
(4.1). The reasoning is exactly as in \cite{P1}, by choosing
$h^2=\frac{e^{-B}}{\lf(4\pi
R^2\ri)^{\frac{n}{2}}}\zeta^2(r_x(y)/R)$, where $\zeta$ be a
nonnegative cut-off function such that $\zeta(t)=1$ for all $t\le
\frac{1}{2}$, and $\zeta(t)=0$ for $t\ge 1$. $B$ is so chosen such
that $\int_M h^2\, dv=1$. Under the assumption (4.3) we have that
$$
\log \frac{V_x(R)}{R^n} +C_1(n)\le B \le \log \frac{V_x(R)}{R^n}
+C_2(n).
$$
Therefore, estimation on the right hand side of (4.2) gives
$$
 -A\le \mu(R^2)\le C_3(n) + B
$$
which implies (4.1) for some $\kappa$. Now argue by contradiction
that (4.1) must holds. If not, we know that (4.3) can not be true.
Namely
$$
V_x(\frac{R}{2}) < \eta V_x(R). \tag 4.4
$$
 We focus on the smaller ball $B_x(\frac{R}{2})$. By the
above argument we would conclude that
$$
V_x(\frac{R}{4}) < \eta V_x(\frac{R}{2}).
$$
Otherwise we would have $V_x(\frac{R}{2})\ge \kappa
\lf(\frac{R}{2}\ri)^n$, which would implies $V_x(\frac{R}{2})\ge
\eta V_x(R)$ by the assumption (4.1) does not hold. Therefore,
iterating the argument we have that
$$
V_x(\frac{R}{2^k})\le \eta^k V_x(R) \tag 4.5
$$
for all natural numbers $k$. This leads to
$$
V_x(r)\le C r^{n\log _23}
$$
for small $r$, which is a contradiction.
\enddemo

The following  result on the diameter of a  manifold with bounded
$\mu(\tau)$ is an easy consequence of Theorem 4.1.

\proclaim{Corollary 4.2} Let $M$ be a Riemanian  manifold such
that $\mu(\tau)\ge -A$ for $1\ge \tau \ge 0$. Assume also that
$V(M)\le V_0$. Then there exists a constant $D=D(A, V_0, n)$ such
that
$$
\text{Diameter}(M) \le D. \tag 4.6
$$
In fact, $D\le 2([\frac{V_0}{\kappa}]+1).$ In particular, it
implies that $M$ is compact if it is complete.
\endproclaim

\proclaim{Concluding remarks} 1) It would be interesting to find
out if there is an interpolation between the entropy formula of
Perelman and (0.4). Namely to find a family monotonicity formulae
connecting both. For the differential Harnack, or Li-Yau-Hamilton
inequality, there is such interpolation in dimension two as shown
by Chow \cite{Ch}. The straightforward formulation seems not to
work. (One could have some differential inequalities connecting
both cases. But the differential inequalities does not give
monotonicity formulae unless on two end points.)

2) It seems that the entropy formula in \cite{P1} is essentially
different from the known one of Hamilton \cite{H1} for the Ricci
flow on Riemann surfaces since it can be used to derive the
uniform scalar curvature bound and diameter bounds without
appealing the Harnack inequality (in \cite{P2} Perelman proved
these results for K\"ahler-Ricci flow with $c_1(M)>0$), unlike the
approach in \cite{H1}, which used the Harnack inequality for the
Ricci flow essentially (in \cite{C-C-Z}, using the similar method
of Hamilton on Riemann surfaces, the authors proved the scalar
curvature and diameter bound for the case when the manifolds has
positive bisectional curvature, which is a special, relatively
easier,  case of what treated in \cite{P2}). Is there any
connection between Perelman's entropy formula and Hamilton's
entropy formula at all?

3) Whether Theorem 0.3 is still true in $n=3$ by assuming instead
the scalar curvature ${\Cal R}(x)$ of $M$ is nonnegative?

4) In \cite{C-N}, the authors proved the matrix Li-Yau-Hamilton
inequality on K\"ahler manifolds with nonnegative bisectional
curvature following an earlier work of Hamilton \cite{H2}, which
can be viewed as a generalized complex Hessian comparison theorem.
The natural question is: does (0.2) have a  matrix version? The
same question applies to Perelman's entropy estimate Corollary
9.3.
\endproclaim


\Refs \widestnumber \key{\bf M-S-Y-1}


\ref\key{\bf B-C-L} \by D. Bakry, D. Concordet and M.  Ledoux
\paper Optimal heat kernel bounds under logarithmic Sobolev
inequalities\jour ESAIM Probab. Statist. \vol 1 \yr 1995/97\pages
391--407
\endref


\ref\key{\bf B} \by W. Beckner \paper Geometric asymptotics and
the logarithmic Sobolev inequality \jour Forum Math.\vol 11\yr1999
\pages 105--137\endref


\ref\key{\bf B-P} \by W. Beckner and M. Pearson \paper On sharp
Sobolev and the logarithmic Sobolev inequality \jour Bull. London.
Math. Soc. \vol 30\yr1998 \pages 80-84\endref


\ref\key{\bf C-N}\by H.D. Cao and L. Ni \paper Matrix
Li-Yau-Hamilton estimates for the heat equation on Kaehler
manifolds \paperinfo Submitted, arXiv: math.DG/0211283
\endref

\ref\key{\bf C-C-Z}\by H.D. Cao,  B-L Chen and X-P Zhu \paper
Ricci flow on compact K\"ahler manifolds of positive bisectional
curvature\paperinfo Preprint
\endref

\ref\key{\bf Ch} \by B. Chow\paper Interpolating between Li-Yau's
and Hamilton's Harnack inequalities on a surface\jour J. Partial
Differential Equations \vol 11 \yr 1998\pages  no. 2, 137--140.
\endref


\ref\key{\bf C-L-Y} \by S.-Y. Cheng, P. Li and S.-T. Yau\paper On
the upper estimate of heat kernel of a complete Riemannian
manifold \jour Amer. J. Math. \vol 103 \yr 1981 \pages 1021--1063
\endref

\ref\key{\bf D} \by E Davies \paper Explicit constants for
Gaussian upper bounds on heat kernels \jour Amer. J. Math. \vol
109 \yr 1987 \pages 319--334
\endref

\ref\key{\bf E} \by K. Ecker\paper Logarithmic Sobolev
inequalities on submanifolds of Euclidean spaces \jour Jour. Reine
Angew. Mat. \vol 552 \yr 2002 \pages 105--118.
\endref


\ref\key{\bf G} \by L. Gross\paper Logarithmic Sobolev
inequalities \jour Amer. J. Math. \vol 97 \yr1975  \pages no. 4,
1061--1083
\endref

\ref\key{\bf H1} \by R. Hamilton \paper The Ricci flow on
surfaces\paperinfo  Mathematics and general relativity (Santa
Cruz, CA, 1986), 237--262, Contemp. Math., ({\bf 71}) , Amer.
Math. Soc., Providence, RI, 1988.
\endref

\ref\key{\bf H2}\by R. Hamilton \paper A matrix Harnack estimate
for the heat equation\jour Comm. Anal. Geom. \vol 1 \yr 1993
\pages 113--126 \endref

\ref\key{\bf Le} \by M. Ledoux \paper On manifolds with
non-negative Ricci curvature and Sobolev inequalities\jour Comm.
Anal. Geom. \vol 7 \yr 1999 \pages 347--353 \endref


\ref\key{\bf L1} \by P. Li \paper On the Sobolev constant and
$p$-spectrum of a compact Riemannian manifold\jour Ann. scient.
\'Ec. Norm. Sup.\vol 13\yr 1980\pages 451--469
\endref

\ref\key{\bf L2} \by P. Li \paper Large time behavior of the heat
equation on complete manifolds with nonnegative Ricci curvature
\jour Ann. of Math. \vol 124 \yr 1986 \pages no. 1, 1--21
\endref

\ref\key{\bf L-Y} \by P. Li and S.-T. Yau\paper On the parabolic
kernel of the Schr\"odinger operator\jour Acta Math.\vol 156\yr
1986\pages 139--168
\endref


\ref\key{\bf N} \by L. Ni\paper Poisson equation and
Hermitian-Einstein metrics on holomorphic vector bundles over
complete noncompact Kahler manifolds \jour Indiana Univ. Math. J.
\vol 51\yr 2002 \pages no 3, 679-704
\endref

\ref\key{\bf P1} \by G. Perelman\paper The entropy formula for the
Ricci flow and its geometric applications \paperinfo arXiv:\
math.DG/\ 0211159
\endref

\ref\key{\bf P2} \by G. Perelman\paper Informal talks and
discussions
\endref

\ref\key{\bf R}\by O. S. Rothaus \paper Analtic inequalities,
isoperimetric inequalities and logarithmic Sobolev inequalities
\jour J. Funct. Anal. \vol 64 \yr 1985 \pages 296--313
\endref


\ref\key{\bf S} \by A. J. Stam \paper Some inequalities satisfied
by the quantities of informations of Fisher and Shannon \jour
Inform. and Control\vol 2\yr 1959 \pages 101--112
\endref

\ref\key{\bf V}\by R. S. Varadhan\paper On the behavior of the
fundamental solution of the heat equation with variable
coefficients \jour Comm. Pure Applied Math. \vol 20 \yr 1967
\pages 431--455
\endref

\ref\key{\bf W}\by F. B. Weissler \paper Logarithmic Sobolev
inequalities for the heat-diffusion semi-group \jour Trans. AMS.
\vol 237\yr 1978 \pages 255--269
\endref



\endRefs
\enddocument